\documentclass[10pt]{article}
\usepackage[dvipdfmx]{graphicx}
\def\qed{\hfill \hbox{${\vcenter{\vbox{
   \hrule height 0.4pt\hbox{\vrule width 0.4pt height 6pt
   \kern5pt\vrule width 0.4pt}\hrule height 0.4pt}}}$}}
\newcommand{\maru}[1]{{\ooalign{\hfil#1\/\hfil\crcr
   \raise.167ex\hbox{\mathhexbox20D}}}}
\newcommand{\bysame}{%
   \leavevmode\hbox to 6em{\hrulefill}\,}

\begin{document}

\centerline{\large\bf On tangle decompositions of twisted torus knots}
\vskip 2mm 

\centerline{by} 
\vskip 2mm 

\centerline{\bf Kanji Morimoto} 
\vskip 5mm 

\centerline{Department of IS and Mathematics, Konan University} 
\vskip -1mm 

\centerline{Okamoto 8-9-1, Higashi-Nada, Kobe 658-8501, Japan}
\vskip -1mm 

\centerline{e-mail : morimoto@konan-u.ac.jp} 
\vskip 5mm 

{\bf Abstract.} \ In the present paper, we will show that for any integer $n > 0$ there are infinitely many twisted torus knots with $n$-string essential tangle decompositions. 
\vskip 5mm 

Keywords and phrases \ : \ twisted torus knots, tangle decompositions 
 
2010 Mathematics Subject Classification \ : \ 57M25, 57M27 
\vskip 10mm 

{\bf 1. Introduction} 

\quad Let $p, q, r, s$ be integers with $p > r > 1$, $q > 0$, gcd$(p, q)=1$, and let $T(p, q)$ be the torus knot of type $(p, q)$ in $S^3$. For the definition of torus knots $T(p, q)$ we refer to [9]. Add $s$ times full twists on mutually parallel $r$-strands in $T(p, q)$. Then according as [2], we call the knot obtained by this operation a twisted torus knot of type $(p, q ; r, s)$ and denote it by $T(p, q ; r, s)$ as illustrated in Figure 1.  

\begin{figure}[htbp]
\hskip 28mm 
\includegraphics[width=6cm]{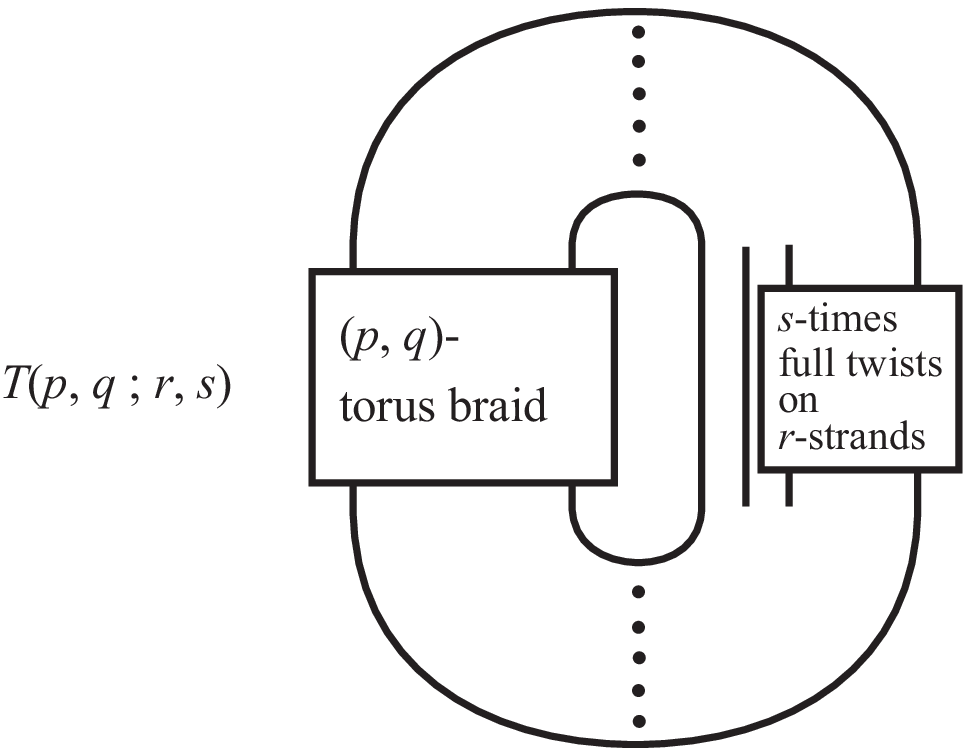}

\center{Figure 1}
\end{figure}

\quad Recently, several interesting results on twisted torus knots have been gotten ([3], [4], [5], [6], [7], [8]). In particular, we have shown in [6] that there are infinitely many composite twisted torus knots. In the present paper, as an extension of the result in [6], we will construct infinitely many twisted torus knots with $n$-string essential tangle decompositions for any integer $n>0$. In fact we will show :  
\vskip 3mm 

{\bf Theorem 1.1.} \ \it Let $e>0, \ k_1>1, \ k_2>1, \ x_1>0, \ x_2>0$ be integers with {\rm gcd}$(x_1, x_2)=1$. Put 

\quad $p=((e+1)(k_1+k_2-1)+1)x_1+(e+1)x_2$, 

\quad $q=(e(k_1+k_2-1)+1)x_1+ex_2$, 

\quad $r=((e+1)(k_1+k_2-1)-k_1+2)x_1+ex_2$ and 

\quad $s=-1$. 

Then we have $:$ 

$(1)$ $T(p, q ; r, s)$ has an $x_1$-string essential tangle decomposition. 

$(2)$ The decomposition is obtained by the $x_1$-string fusion of the torus knot $T((k_1-1)x_1+x_2, e((k_1-1)x_1+x_2)+x_1)$ and the torus link $T(k_2x_1, -((e+1)k_2+1)x_1)$. 

$(3)$ $T(p, q ; r, s)$ has an essential torus in the exterior whose companion is the torus knot  $T(k_2, -(e+1)k_2-1)$. 

Therefore, for any integer $n>0$, by putting $x_1=n$ we get infinitely many twisted torus knots with $n$-string essential tangle decompositions. 
\rm 
\vskip 3mm 

{\bf Example 1.2.} \ Put $e=1, k_1=k_2=2, x_1=2, x_2=3$. Then by Theorem 1.1, we see that $T(20, 11 ; 15, -1)$ has a 2-string essential tangle decomposition which is obtained by the 2-string fusion of $T(5, 7)$ and $T(4, -10)$ as in Figure 2 (c.f. Example 3.3). By tubing the decomposing 2-sphere along the torus link $T(4, -10)$, we have an essential torus whose companion is the torus knot $T(2, -5)$. 

\begin{figure}[htbp]
\hskip 15mm 
\includegraphics[width=9.5cm]{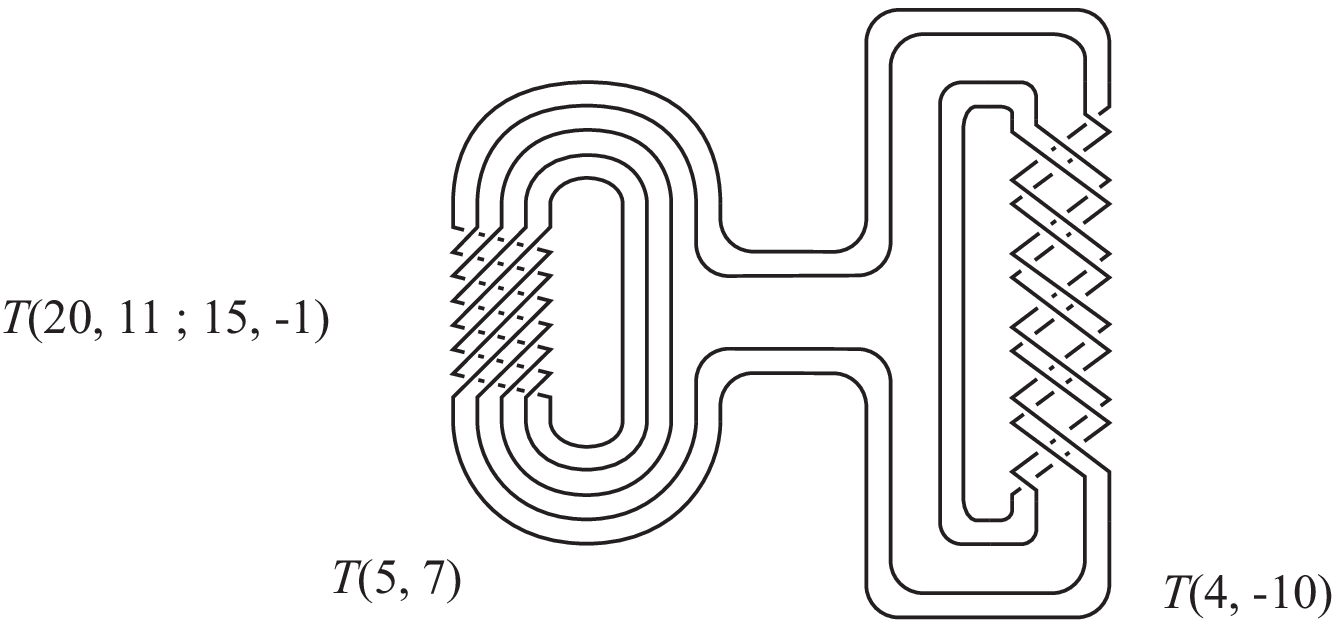}

\center{Figure 2}
\end{figure}

{\bf Remark 1.3.} \ Suppose $x_1=1$ in Theorem 1.1. Then by putting $k_1'=k_1+x_2-1$ and $k_2'=k_2$, we have :

\quad $p=(e+1)(k_1+k_2-1)+1+(e+1)x_2=(e+1)(k_1'+k_2')+1$, 

\quad $q=e(k_1+k_2-1)+1+ex_2=e(k_1'+k_2')+1$, 

\quad $r=(e+1)(k_1+k_2-1)-k_1+2+ex_2=p_0-k_1'$ and 

\quad $s=-1$. 

Then by Theorem 1 of [6], we see that $T(p, q ; r, s)$ is a composite twisted torus knot if $x_1=1$. See Theorem 3.1 in the present paper.  
\vskip 3mm  

\quad In [4], S. Y. Lee showed the following (c.f. [8]) :  
\vskip 3mm 

{\bf Theorem 1.4 (Theorem 1 of [4]).} \ \it Suppose $r=qk$ for some integer $k$. Then  $T(p, q ; r, s)$ is the $(q, p+k^2qs)$-cable knot along the torus knot $T(k, ks+1)$. \rm 
\vskip 3mm 

\quad Thus, by noting (3) of Theorem 1.1 and Remark 1.3, we conjecture the following : 
\vskip 3mm 

{\bf Conjecture.} \ \it A twisted torus knot $T(p, q ; r, s)$ has an essential torus in the exterior, if and olny if it is a knot in Theorem $1.1$ or $r=qk$ for some integer $k$. In particular, $T(p, q ; r, s)$ is a composite knot if and only if it is a knot in Theorem $1.1$ with $x_1=1$. \rm 
\vskip 3mm 

\quad Throughout the present paper, we will work in the piecewise linear category. For a manifold $X$ and a subcomplex $Y$ in $X$, we denote a regular neighborhood of $Y$ in $X$ by $N(Y, X)$ or $N(Y)$ simply. 
\vskip 3mm 

{\bf 2. Parallelized torus knots and parallelized twisted torus knots} 

\quad Let $T(p_0, q_0)$ be the torus knot of type $(p_0, q_0)$, where $p_0$ and $q_0$ are positive coprime integers with $p_0>1$, and let $x_1$ and $x_2$ be positive integers. Take four points P$_1$, P$_2$, P$_3$ and P$_4$ on the adjacent two strands in $T(p_0, q_0)$ as in Figure 3. Then replace the arc P$_1$ through P$_3$ with $x_1$ parallel strings and the arc P$_2$ through P$_4$ with $x_2$ parallel strings. In addition, replace the rectangle P$_1$P$_2$P$_3$P$_4$ with $x_1+x_2$ strands as in Figure 4. Then we get a torus knot or a torus link $T(p, q)$ for some $p, q$. 

\begin{figure}[htbp]
\hskip 30mm 
\includegraphics[width=5.5cm]{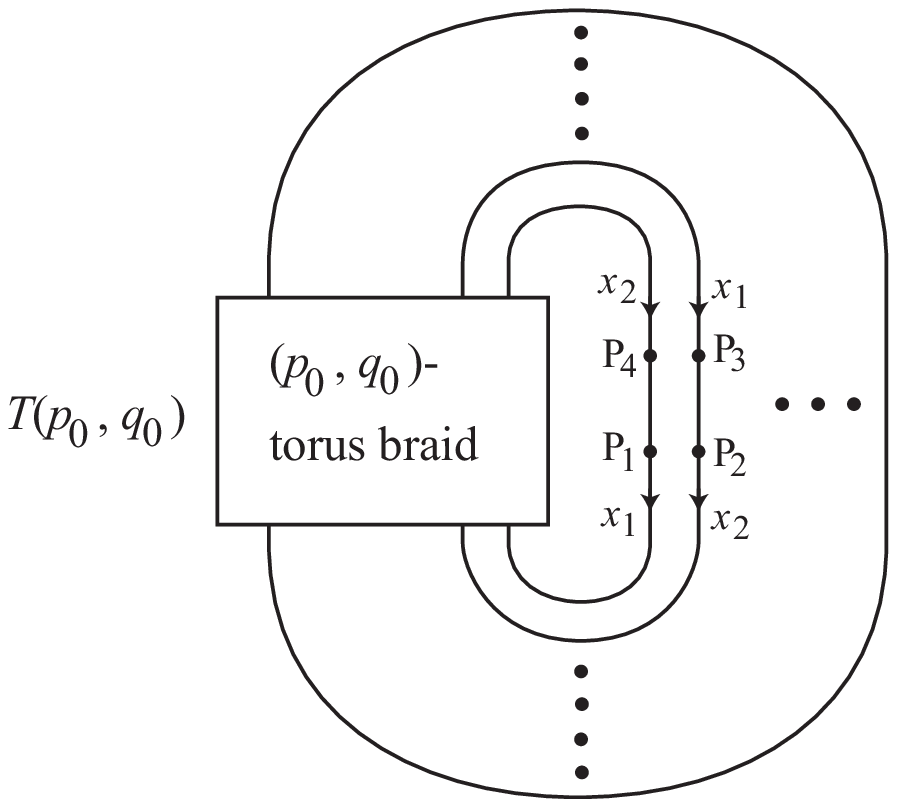}

\center{Figure 3}
\end{figure}

\begin{figure}[htbp]
\vskip 5mm 
\hskip 30mm 
\includegraphics[width=7cm]{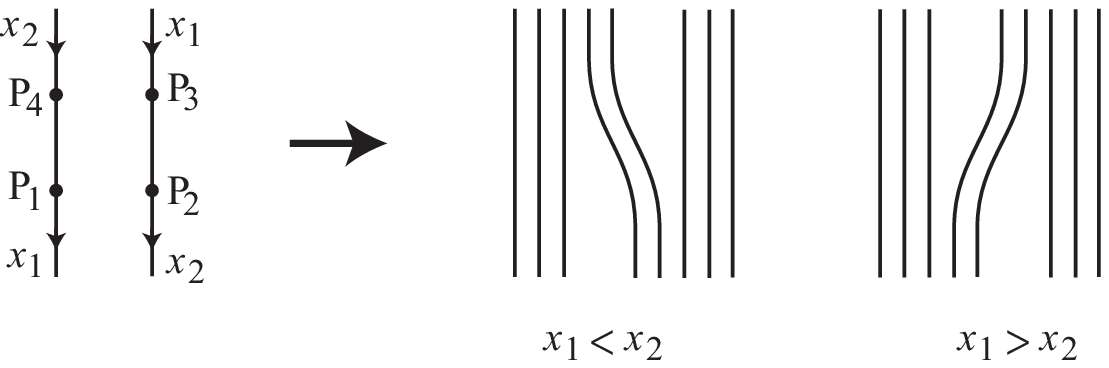}

\center{Figure 4}
\end{figure}

\begin{figure}[htbp]
\hskip 35mm 
\includegraphics[width=6cm]{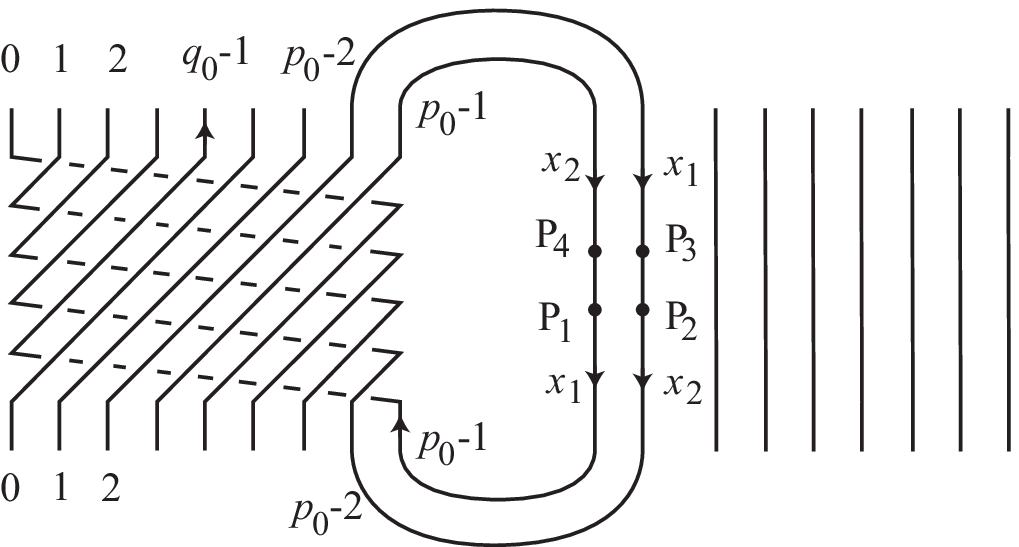}

\center{Figure 5}
\end{figure}

\quad Let's detect $p$ and $q$. First, number the $p_0$ strings below the $(p_0, q_0)$ torus braid $0, 1, 2, \cdots, p_0-2, p_0-1$ as in Figure 5. The arc starting at P$_1$ goes into the braid at $p_0-1$ and goes out at $q_0-1$. After coming back round once it goes into the braid again and goes out at $2q_0-1$. Next it goes out the braid at $3q_0-1$. By continuing these procedures, it finally goes out at $aq_0-1 \equiv p_0-2$ (mod $p_0$). Then it meets the point P$_3$. Hence we have $aq_0 \equiv -1$ (mod $p_0$). Similarly the arc starting at P$_2$ goes into the braid at $p_0-2$ and goes out at $q_0-2$. Then it goes out the braid at $2q_0-2, 3q_0-2, \cdots $, and finally goes out at $bq_0-2 \equiv p_0-1$ (mod $p_0$). Then it meets the point P$_4$. Hence we have $bq_0 \equiv 1$ (mod $p_0$). 

\quad Thus we have $p=ax_1+bx_2$, where $a$ and $b$ are the least positive integers such that  $aq_0 \equiv -1$ (mod $p_0$), $bq_0 \equiv 1$ (mod $p_0$) and $a+b=p_0$. 

\quad By the similar arguments, we have  $q=cx_1+dx_2$, where $c$ and $d$ are the least positive integers such that  $cp_0 \equiv 1$ (mod $q_0$), $dp_0 \equiv -1$ (mod $q_0$) and $c+d=q_0$. 

\quad In general, we have :  
\vskip 3mm 

{\bf Proposition 2.1.} \ \it For coprime positive integers $p_0$ and $q_0$, there uniquely exist positive integers $a, b, c, d$ which satisfy the following conditions $:$ 
\vskip 3mm
 
\centerline{
$(1)$ $\cases{a+b=p_0 \cr aq_0 \equiv -1 \ ({\rm mod} \ p_0) \cr bq_0 \equiv \ \ 1 \ ({\rm mod} \ p_0) }$ 
\hskip 10mm 
$(2)$ $\cases{c+d=q \cr cp_0 \equiv \ \ 1 \ ({\rm mod} \ q_0) \cr dp_0 \equiv -1 \ ({\rm mod} \ q_0) }$} \rm 
\vskip 3mm 

{\bf Proof.} \ Consider the set $\{ 0, q_0, 2q_0, \cdots, (p_0-1)q_0 \}$. Then, by gcd$(p_0, q_0)=1$, these $p_0$ integers are different to each other (mod $p_0$). Then this set coincides with the set $\{ 0, 1, 2, \cdots, p_0-1 \}$ (mod $p_0$), and hence there is only one integer $a$ with $aq_0 \equiv p_0-1 \equiv -1$  (mod $p_0$). Then put $b=p_0-a$, and we have $bq_0=(p_0-a)q_0=p_0q_0-aq_0 \equiv 0-(-1)=1$ (mod $p_0$). This completes the proof of (1), and (2) is proved similarly. \qed 
\vskip 3mm 

\quad Under the above situations, we have : 
\vskip 3mm 

{\bf Proposition 2.2.} \ \it Let $x_1$ and $x_2$ be positive integers, and put  $p=ax_1+bx_2$ and  $q=cx_1+dx_2$. Then ${\rm gcd}(p, q) = {\rm gcd}(x_1, x_2)$. In particular, $T(p, q)$ is a torus knot if and only if ${\rm gcd}(x_1, x_2)=1$. \rm 

{\bf Proof.} \ Put gcd$(x_1, x_2)=k$. Then we can put $x_1=ky_1, \ x_2=ky_2$, and $p=k(ay_1+by_2), \ q=k(cy_1+dy_2)$. Hence gcd$(p, q) \ge k = $ gcd$(x_1, x_2)$. 

\quad Conversely, put gcd$(p, q)=k$. Then we can put $p=kp_1, \ q=kq_1$. 
\vskip 3mm 

\quad Since \ 
$\left[ \begin{array}{c} p \\ q \end{array} \right] = \left[ \begin{array}{cc} a & b \\ c & d \end{array} \right] \left[ \begin{array}{c} x_1 \\ x_2 \end{array} \right]$, \ we have :   
\vskip 3mm 

\centerline{$\displaystyle \left[ \begin{array}{c} x_1 \\ x_2 \end{array} \right] = {1 \over {ad-bc}}\left[ \begin{array}{rr} d & -b \\ -c & a \end{array} \right] \left[ \begin{array}{c} kp_1 \\ kq_1 \end{array} \right] \cdots \maru{1}$}   
\vskip 3mm 

\quad Then $|ad-bc|<p_0q_0-1$ because $0<a, \ b<p_0$ and $0<c, \ d<q_0$. Moreover $ad-bc=a(q_0-c)-(p_0-a)c=aq_0-cp_0 \equiv -1$ (mod $p_0$), (mod $q_0$). This implies that $ad-bc=-1$, and by \maru{1} we have gcd$(x_1, x_2) \ge k =$ gcd$(p, q)$. This completes the proof. \qed 
\vskip 3mm 

\quad Let $T(p_0, q_0 ; r_0, s_0)$ be a given twisted torus knot. Then by the same way as the case of torus knots, we can construct a parallelized twisted torus knot or a parallelized twisted torus link $T(p, q ; r, s)$. Then $p=ax_1+bx_2, \ q=cx_1+dx_2, \ r=r_1x_1+r_2x_2$ and $s=s_0$, where $a, b, c, d$ are those integers in Proposition 2.1 and $r_1, r_2$ are some positive integers with $r_1+r_2=r_0$. Then by Proposition 2.2, we see that  $T(p, q ; r, s)$ is a knot if and only if gcd$(x_1, x_2)=1$. 
\vskip 3mm 

{\bf 3. Proof of Theorem 1.1} 

\quad Let $B$ be a 3-ball and $t^1 \cup t^2 \cup \cdots \cup t^n$ $n$ arcs properly embedded in $B$, then we call the pair $(B, t^1 \cup t^2 \cup \cdots \cup t^n)$ an $n$-string tangle. We say that $(B, t^1 \cup t^2 \cup \cdots \cup t^n)$ is essential if $cl(\partial B - N(t^1 \cup t^2 \cup \cdots \cup t^n))$ is incompressible in $cl(B - N(t^1 \cup t^2 \cup \cdots \cup t^n))$ if $n > 1$, and $t^1$ is a knotted arc in $B$ if $n = 1$, where $N(t^1 \cup t^2 \cup \cdots \cup t^n)$ is a regular  neighborhood of $t^1 \cup t^2 \cup \cdots \cup t^n$ in $B$, and that the tangle is inessential if it is not essential. We say that a knot $K$ in the 3-sphere $S^3$ has an $n$-string essential tangle decomposition if $(S^3, K)$ is decomposed into two $n$-string essential tangles $(B_1, t_1^1 \cup t_1^2 \cup \cdots \cup t_1^n) \cup (B_2, t_2^1 \cup t_2^2 \cup \cdots \cup t_2^n)$, and that the decomposition is inessential if it is not essential. 

\quad To prove Theorem 1.1, we construct parallelized twisted torus knots from composite twisted torus knots, and we will show that the decomposing 2-sphere of the connected sum becomes the decomposing 2-sphere of the tangle decomposition. In [6], we have shown the following: 
\vskip 3mm 

{\bf Theorem 3.1 (Theorem 1 of [6]).} \it Let $e>0, k_1>1, k_2>1$ be integers, and put $p_0=(e+1)(k_1+k_2)+1$, $q_0=e(k_1+k_2)+1$, $r_0=p_0-k_1$ and $s_0=-1$. Then $T(p_0, q_0 ; r_0, s_0)$ is the connected sum of $T(k_1, ek_1+1)$ and $T(k_2, -(e+1)k_2-1)$. \rm 
\vskip 3mm 

\quad To get parallelized twisted torus knots from the above composite knots, first we calculate the integers $a, b, c, d$ in Proposition 2.1 to get $p$ and $q$. 
\vskip 3mm 

{\bf Proposition 3.2.} \ \it For $p_0=(e+1)(k_1+k_2)+1$ and $q_0=e(k_1+k_2)+1$, those integers $a, b, c, d$ in Proposition $2.1$ are as follows $:$ 

$a=(e+1)(k_1+k_2-1)+1$, $b=e+1$, $c=e(k_1+k_2-1)+1$ and $d=e$. \rm 

{\bf Proof.} \ First we have $b=e+1$, because $(e+1)q_0=(e+1)(e(k_1+k_2)+1)=(e+1)e(k_1+k_2)+e+1=e((e+1)(k_1+k_2)+1)+1=ep_0+1 \equiv 1$  (mod $p_0$). Then $a=p_0-b=(e+1)(k_1+k_2)+1-(e+1)=(e+1)(k_1+k_2-1)+1$ and $aq_0=(p_0-b)q_0=p_0q_0-bq_0 \equiv -1$ (mod $p_0$). 

\quad Next we have $d=e$ because $ep_0=e((e+1)(k_1+k_2)+1)=e(e+1)(k_1+k_2)+e=(e+1)e(k_1+k_2)+e+1-1=(e+1)(e(k_1+k_2)+1)-1=(e+1)q-1 \equiv -1$ (mod $q_0$). Then $c=q_0-d=e(k_1+k_2)+1-e=e(k_1+k_2-1)+1$ and $cp_0=(q_0-d)p_0=q_0p_0-dp_0  \equiv 1$ (mod $q_0$). This completes the proof. \qed 
\vskip 3mm 

\quad To calculate $r$ and to get concrete expression of the tangle decompositions, consider the twisted torus knots in Theorem 3.1. Put $K_0=T(p_0, q_0 ; r_0, s_0), \ K_1=T(k_1, ek_1+1)$ and $K_2=T(k_2, -(e+1)k_2-1)$, then $K_0=K_1 \# K_2$. Let $V$ be a standard genus two handlebody in $S^3$, and put $\partial V = F$. Then, since any twisted torus knot can be embedded in $F$ in a standard way, we may assume that $K_0$ is in $F$. Let $S$ be the decomposing 2-sphere of the connected sum $K_0 = K_1 \# K_2$, then, by the proof in [6] of Theorem 3.1, we may assume that $S$ intersects $V$ in a separating disk and that $S \cap F$ is a single loop. Then, by noting that $p_0-r_0=k_1$, $S \cap F$ runs along the both sides of $k_1$ strings and $(S \cap F) \cap K_0$ consists of the two points Q$_1$ and Q$_2$ as in Figure 6, where Figure 6 is the case of $e=1, k_1=3, k_2=2$ and the connected sum is $T(11, 6 ; 8, -1) = T(3, 4) \# T(2, -5)$.    

\begin{figure}[htbp]
\hskip 35mm 
\includegraphics[width=7cm]{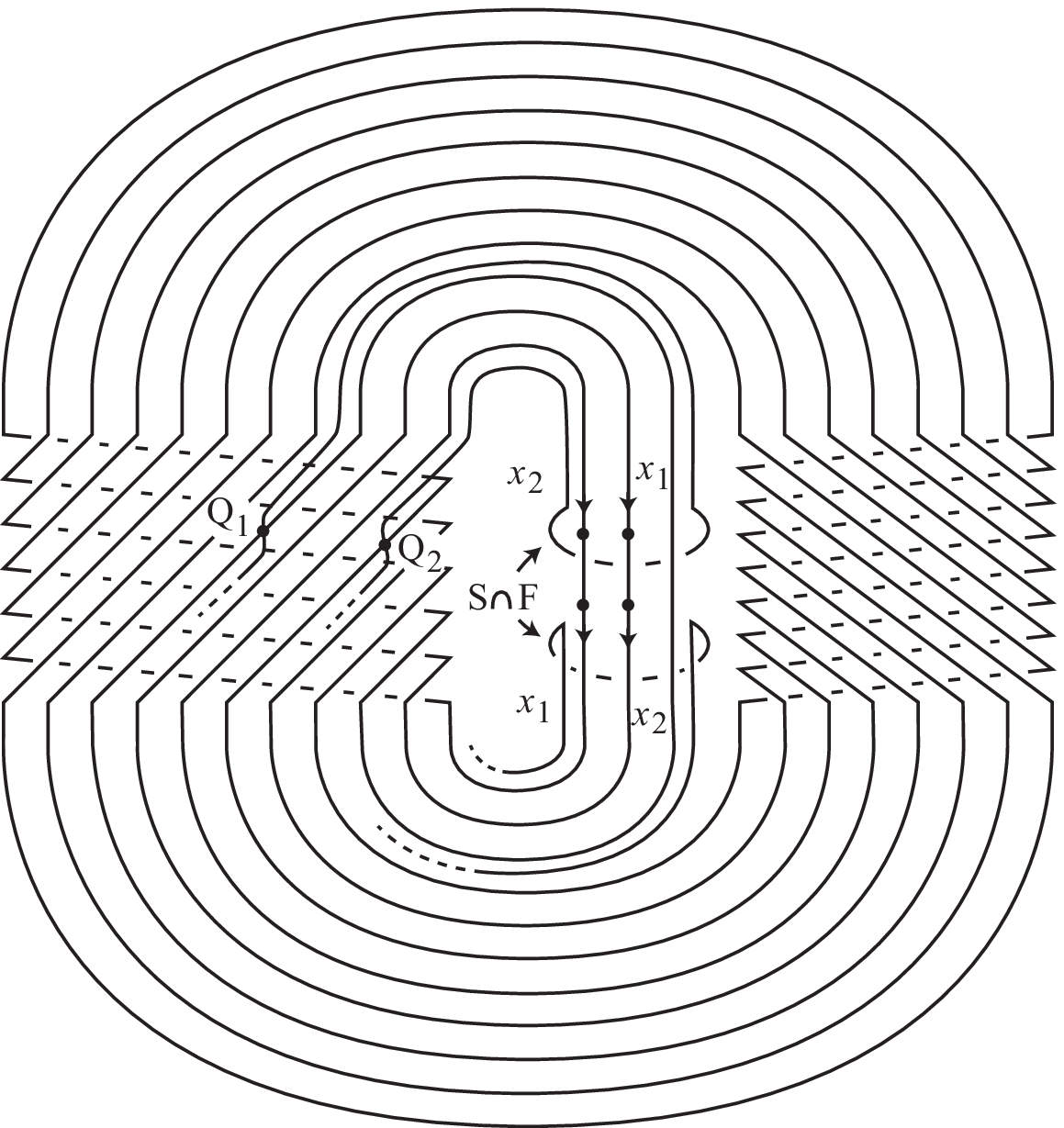}

\center{Figure 6}
\end{figure}

\quad We split $K_0$ at  Q$_1$, Q$_2$ into two arcs, and to get $K_1$ and $K_2$, connect the two points with the arc in the disk $S \cap V$. First we consider $K_1$ as in Figure 7 $(k_1=5)$. Then, by noting that $p_0-r_0=k_1$ and $K_1=T(k_1, ek_1+1)$, we see that the arc P$_2$P$_4$ is contained in $K_1$, and hence exactly one string of $k_1$ strings is replaced with $x_2$ parallel strings. Then the other $(k_1-1)$ strings are contained in the arc P$_1$P$_3$ and are replaced with $x_1$ parallel strings. Thus we get the torus knot $T((k_1-1)x_1+x_2, e((k_1-1)x_1+x_2)+x_1)$, and this torus knot intersects the original decomposing 2-sphere in $x_1$ points at each Q$_i (i=1,2)$. This implies that $r=p-((k_1-1)x_1+x_2)$. 

\begin{figure}[htbp]
\hskip 45mm 
\includegraphics[width=4cm]{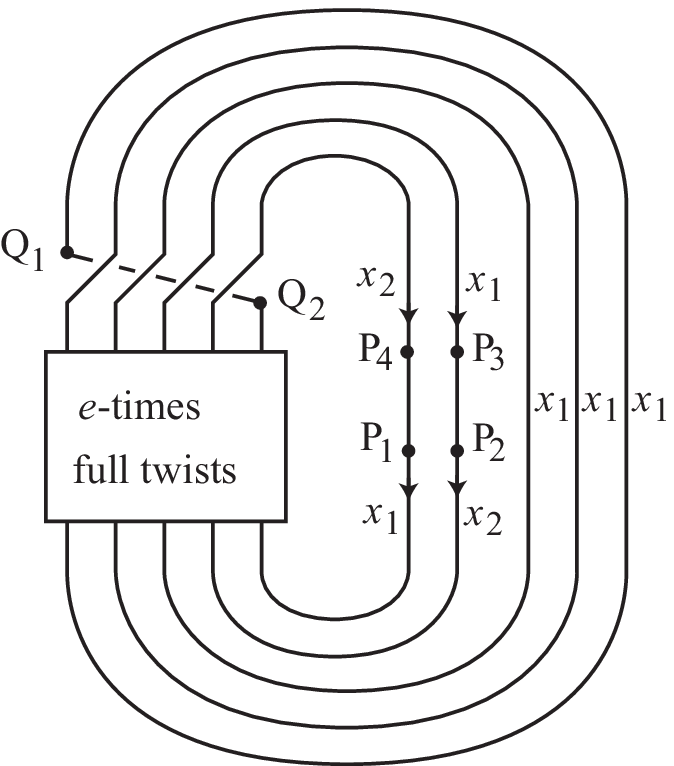}

\center{Figure 7}
\end{figure}

\quad For the knot $K_2$, by the above arguments, we see that the whole string of $K_2$ is contained in the arc  P$_1$P$_3$. Hence by replacing the whole string with $x_1$ parallel strings, we get the torus link 
 $T(k_2x_1, -((e+1)k_2+1)x_1)$. This torus link intersects the original decomposing 2-sphere in $x_1$ points at each  Q$_i (i=1,2)$ similarly to the case of $K_1$. 

\quad By summarizing the above calculations, we have the following, and get the knots in Theorem 1.1.  

\quad $p=ax_1+bx_2=((e+1)(k_1+k_2-1)+1)x_1+(e+1)x_2$

\quad $q=cx_1+dx_2=(e(k_1+k_2-1)+1)x_1+ex_2$

\quad $r=p-((k_1-1)x_1+x_2)=((e+1)(k_1+k_2-1)+1)x_1+(e+1)x_2-((k_1-1)x_1+x_2)$ 

\hskip 8mm $=((e+1)(k_1+k_2-1)-k_1+2)x_1+ex_2$

\quad $s=s_0=-1$ 

\quad Finally, we need to show that the above tangle decompositions are all essential. If $x_1=1$, then the decompositions are the connected sums and are all essential because both $k_1$ and $k_2$ are greater than one and factor knots are non-trivial knots.  

\quad Suppose $x_1>1$. By the definition of tangles, we see that an $n$-string tangle $(B, t^1 \cup t^2 \cup \cdots \cup t^n)$ with $n>1$ is essential if and only if there is no disk properly embedded in $B$ which separates the arcs $t^1 \cup t^2 \cup \cdots \cup t^n$. From this view point, in the next section, we will show that both of $x_1$-string tangles constructed from the torus knot $T((k_1-1)x_1+x_2, e((k_1-1)x_1+x_2)+x_1)$ and the torus link $T(k_2x_1, -((e+1)k_2+1)x_1)$ are essential (Proposition 4.3). 

\quad In addition, by tubing the decomposing 2-sphere along the torus link $T(k_2x_1, -((e+1)k_2+1)x_1)$, we have an essential torus whose companion is the torus knot $T(k_2, -(e+1)k_2-1)$. This completes the proof of Theorem 1.1. \qed 
\vskip 3mm 

{\bf Example 3.3.} \ Put $e=1, k_1=k_2=2$, and let $x_1, x_2$ be positive integers. Then by the above arguments, $p=7x_1+2x_2, q=4x_1+x_2, r=6x_1+x_2$ and $T(p, q ; r, -1)$ is the $x_1$-string fusion of $T(x_1+x_2, 2x_1+x_2)$ and $T(2x_1, -5x_1)$. Hence by putting $x_1=2, x_2=3$, we see that  $T(20, 11 ; 15, -1)$ is the 2-string fusion of $T(5, 7)$ and $T(4, -10)$, and this is the example in Introduction. We note that there may be ambiguity on twists in the fusion. But by more detailed arguments we get Figure 2.  If we put $x_1=2$ and $x_2=1$, then we see that $T(16, 9 ; 13, -1)$ is the 2-string fusion of $T(3, 5)$ and $T(4, -10)$. This is the smallest example of our knots. 
\vskip 3mm 

{\bf 4. Essential tangles} 

\quad Let $p, q$ be coprime integers with $2<p<q$, and $k$ an integer with $1<k<p$. Consider the torus knot $T(p, q)$ and take an arc $\alpha$ which intersects $k$ strings in the parallel $p$ strings of $T(p, q)$ as in Figure 8(1). Let $N(\alpha)$ be a regular neighborhood of $\alpha$ in $S^3$. Put $B=cl(S^3-N(\alpha))$ and $t(p, q; k)=cl(T(p, q)-N(\alpha))$. Then the pair $(B, t(p, q; k))$ is a $k$-string tangle as in Figure 8(2). 

\begin{figure}[htbp]
\hskip 25mm 
\includegraphics[width=8cm]{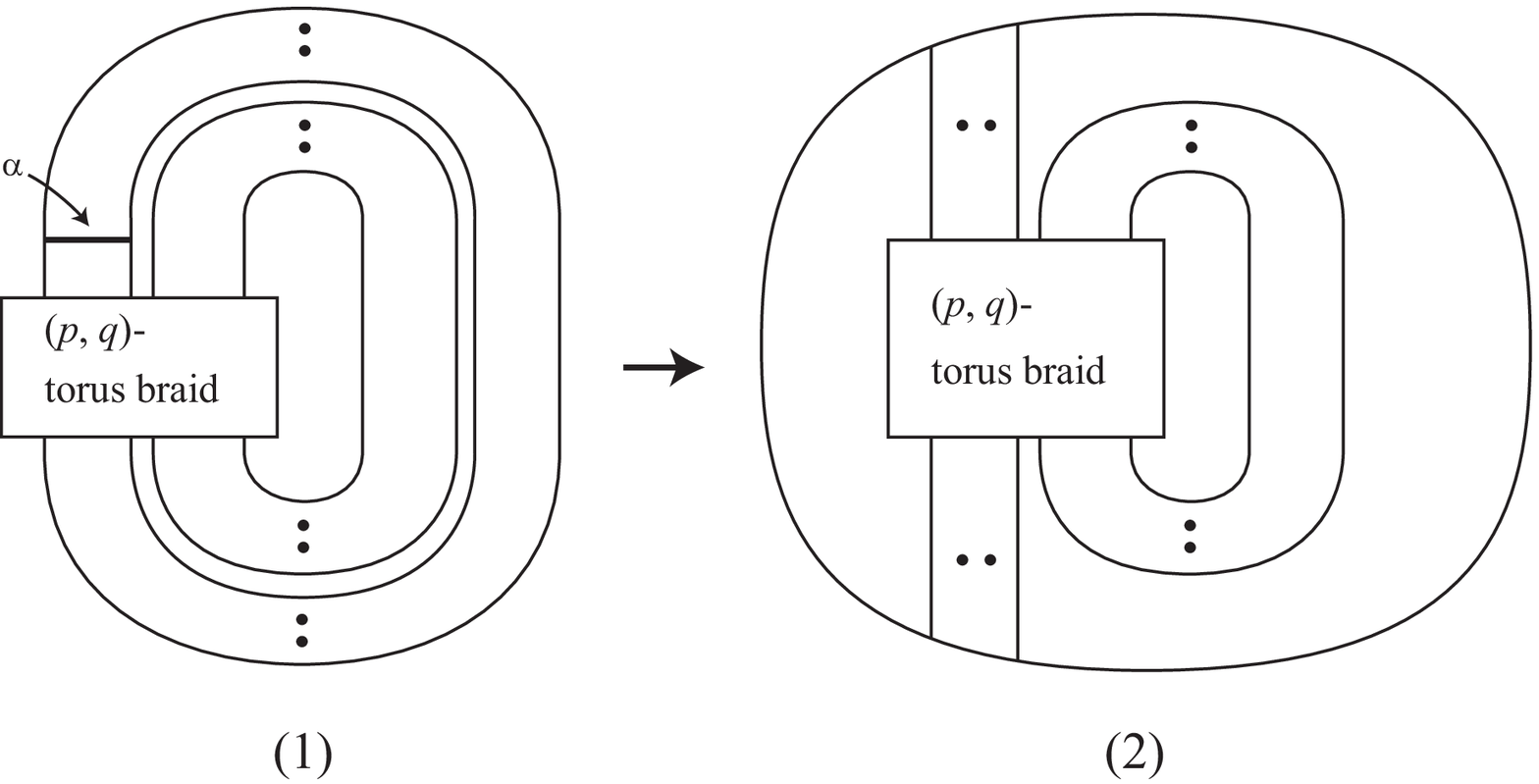}

\center{Figure 8}
\end{figure}

\vskip 3mm 

{\bf Lemma 4.1.} \ \it If $k=2$, then the tangle $(B, t(p, q ; 2))$ is an essential tangle. \rm 

{\bf Proof.} \ Put $t(p, q ; 2)=t_1 \cup t_2$. We first show that at least one of $t_1$ and $t_2$ is a knotted arc in $B$. By $p<q$, we can put $q=np+m \ (0<m<p)$. Then, since $t_1$ and $t_2$ are arcs properly embedded in $B$ each of which is a local torus knot, we can put $t_1=T(a, na+c)$ and $t_2=T(b, nb+d)$, where $a+b=p$ and $c+d=m$ by the arguments similar to those in Section 2. We note that $(na+c)+(nb+d)=n(a+b)+(c+d)=np+m=q$. Then, by $a \le na+c$, $b \le nb+d$ and $a+b=p>2$, we see that at least one of $t_1$ and $t_2$ is a knotted arc in $B$. 

\quad Suppose the tangle $(B, t(p, q ; 2))$ is inessential. Then, by the definition of essential tangles, there is a disk properly embedded in $B$ which separates $t_1$ and $t_2$. Hence $cl(B-N(t_1 \cup t_2))$ is not a handlebody. However, by [1], the arc $\alpha$ connecting the adjacent strings is an unknotting tunnel of $T(p, q)$, and hence $cl(B-N(t_1 \cup t_2))$ is a handlebody. This is a contradiction, and completes the proof. \qed 
\vskip 3mm 

{\bf Lemma 4.2.} \ \it $(B, t(p, q ; k))$ has exactly two parallel classes of strings. \rm 

{\bf Proof.} \ First suppose $k=3$, and number the 3 strings under the torus braid 0, 1 and 2. Then the 3 strings above the torus braid has the three cases, (0, 1, 2), (2, 0, 1) and (1, 2, 0). However if it is (0, 1, 2), then $T(p, q)$ is a 3-component link. Hence we have the two cases (2, 0, 1) and (1, 2, 0) as in Figure 9. Then we can see that the adjacent two strings (0, 1) or (1, 2) are mutually parallel strings.   
 
\begin{figure}[htbp]
\hskip 40mm 
\includegraphics[width=5.5cm]{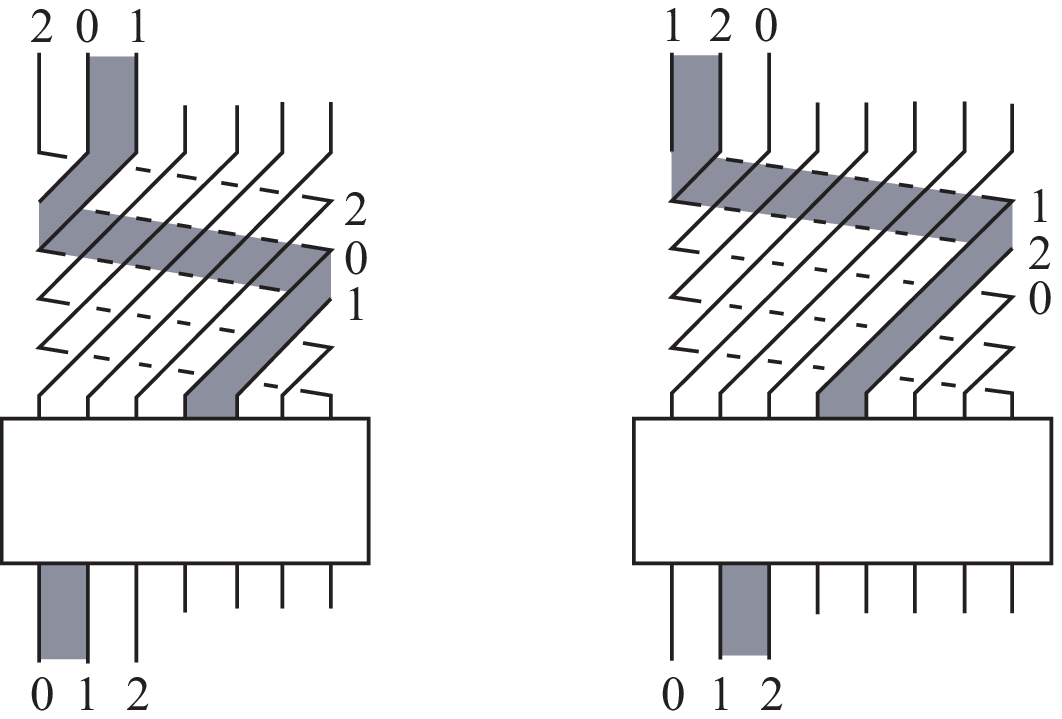}

\center{Figure 9}
\end{figure}

\quad Next suppose $k>3$. Then by the arguments similar to the case of $k=3$, we see that there are adjacent strings which are mutually parallel. Then by inductive arguments, we see that the parallel classes of the strings consists of two classes. \qed 
\vskip 3mm 

{\bf Proposition 4.3.} \ \it The tangle $(B, t(p, q ; k))$ is an essential tangle. \rm 

{\bf Proof.} \ By Lemma 4.2, $t(p, q ; k)=t_1 \cup t_2 \cup \cdots \cup t_k$ has two parallel classes. Hence we may assume that there are two rectangles $R_1, R_2$ in $B$ such that $t_1 \cup \cdots \cup t_i \subset R_1$ and $t_{i+1} \cup \cdots \cup t_k \subset R_2$, where $\partial R_1$ consists of $t_1 \cup t_i \cup ({\rm two \ arcs \ in} \ \partial B)$ and $\partial R_2$ consists of $t_{i+1} \cup t_k \cup ({\rm two \ arcs \ in} \ \partial B)$. 

\quad Suppose the tangle $(B, t(p, q ; k))$ is inessential. Then, by the definition of essential tangles, there is a disk, say $D$, properly embedded in $B$ which separates $t_1 \cup \cdots \cup t_k$. Suppose $R_1 \cap D \ne \emptyset$. Then by standard cut and paste operations and by replacing disks if necessary, we may assume that each component of $R_1 \cap D$ is an arc parallel to the strings $t_1 \cup t_i$. This means that those strings are all unknotted. Then by some cut and paste operations, we can get another separating disk, say $D'$, with $D' \cap R_1 = \emptyset$. If $R_2 \cap D \ne \emptyset$, we can remove the intersections similarly. 

\quad Hence we may assume that $(R_1 \cup R_2) \cap D = \emptyset$ and that $D$ separates the two rectangles $R_1$ and $R_2$. Then by the arguments similar to those in the proof of Lemma 4.1, we have a contradiction. This completes the proof. \qed

\vskip 5mm 

{\bf References} 

[1] \ M. Boileau, M. Rost and H. Zieschang, \ {\it On Heegaard decompositions of torus} 
\vskip -1mm \hskip 6mm {\it knot exteriors and related Seifert fibre spaces}, Math. Ann. {\bf 279} (1988) 553-581.
\vskip -1mm 
[2] \ P. J. Callahan, J. C. Dean and J. R.Weeks,
{\it The simplest hyperbolic knots}, 
\vskip -1mm \hskip 6mm 
J. Knot Theory Ramifications, {\bf 8} (1999) 279--297. 
\vskip -1mm 
[3] \ J. H. Lee, \ {\it Twisted torus knots $T(p, q ; 3, s)$ are tunnel number one},  
\vskip -1mm \hskip 6mm 
J. Knot Theory Ramifications, {\bf 20} (2011) 807--811.
\vskip -1mm 
[4] \ S. Y. Lee, \ {\it Twisted torus knots $T(p, q ; kq, s)$ are cable knots},  
\vskip -1mm \hskip 6mm 
J. Knot Theory Ramifications, {\bf 21} (2012) 1250005(1-4).
\vskip -1mm 
[5] \ K. Morimoto, \ {\it Essential surfaces in the exteriors of torus knots with twists on} 
\vskip -1mm \hskip 6mm 
{\it $2$-strands}, preprint. 
\vskip -1mm 
[6] \ K. Morimoto, \ {\it On composite twisted torus knots}, to appear in Tokyo Journal 
\vskip -1mm \hskip 6mm 
of Mathematics. 
\vskip -1mm 
[7] \ K. Morimoto, M. Sakuma and Y. Yokota, \ {\it Examples of tunnel number one} 
\vskip -1mm \hskip 6mm
{\it knots which have the property `` $1 + 1 = 3$ ''}, \ Math. Proc. Camb. Phil. Soc. 
\vskip -1mm \hskip 6mm 
{\bf 119} (1996) 113--118. 
\vskip -1mm 
[8] \ K. Morimoto and Y. Yamada, \ {\it A note on essential tori in the exteriors of torus} 
\vskip -1mm \hskip 6mm 
{\it knots with twists}, Kobe J. Math., {\bf 26} (2009) 29--34. 
\vskip -1mm 
[9] \ D. Rolfsen, \ {\it Knots and Links}, AMS Chelsea Publishing (2003).

\end{document}